\documentclass[11pt]{amsart}
\input{diagrams}

\newtheorem{proposition}{Proposition}[section]
\newtheorem{lemma}[proposition]{Lemma}
\newtheorem{corollary}[proposition]{Corollary}
\newtheorem{theorem}[proposition]{Theorem}

\theoremstyle{definition}

\newtheorem{example}[proposition]{Example}

\theoremstyle{remark}
\newtheorem{remark}[proposition]{Remark}

\newcommand{\thlabel}[1]{\label{th:#1}}

\newcommand{\selabel}[1]{\label{se:#1}}
\newcommand{\seref}[1]{Section~\ref{se:#1}}
\newcommand{\lelabel}[1]{\label{le:#1}}
\newcommand{\leref}[1]{Lemma~\ref{le:#1}}
\newcommand{\prlabel}[1]{\label{pr:#1}}
\newcommand{\prref}[1]{Proposition~\ref{pr:#1}}
\newcommand{\colabel}[1]{\label{co:#1}}

\newcommand{\relabel}[1]{\label{re:#1}}
\newcommand{\reref}[1]{Remark~\ref{re:#1}}
\newcommand{\exlabel}[1]{\label{ex:#1}}

\newcommand{\eqlabel}[1]{\label{eq:#1}}
\newcommand{\equref}[1]{(\ref{eq:#1})}

\newcommand{\Pic}{{\rm Pic}}
\newcommand{\Hom}{{\rm Hom}}

\newcommand{\Aut}{{\rm Aut}\,}

\newcommand{\Ker}{{\rm Ker}\,}

\newcommand{\im}{{\rm Im}\,}

\newcommand{\Spec}{{\rm Spec}\,}

\def\ot{\otimes}

\def\ZZ{{\mathbb Z}}

\def\GG{{\mathbb G}}

\def\CCC{{\mathfrak C}}

\newcommand{\Cc}{\mathcal{C}}
\newcommand{\Dd}{\mathcal{D}}

\newcommand{\Mm}{\mathcal{M}}

\def\text#1{{\rm {\rm #1}}}

\def\ul{\underline}
\def\dul#1{\underline{\underline{#1}}}

\begin{document}
\title[The relative Picard group]{The relative Picard group of a comodule algebra
and Harrison cohomology}
\author{S. Caenepeel}
\address{Faculty of Engineering Sciences,
Vrije Universiteit Brussel, VUB, B-1050 Brussels, Belgium}
\email{scaenepe@vub.ac.be}
\urladdr{http://homepages.vub.ac.be/\~{}scaenepe/}
\author{T. Gu\'ed\'enon}
\address{Faculty of Engineering Sciences,
Vrije Universiteit Brussel, VUB, B-1050 Brussels, Belgium}
\email{tguedeno@vub.ac.be}
\thanks{Research supported by the project G.0278.01 ``Construction
and applications of non-commutative geometry: from algebra to physics"
from FWO Vlaanderen}
\subjclass{16W30}
\keywords{}
\begin{abstract}
Let $A$ be a commutative comodule algebra over a commutative bialgebra $H$.
The group of invertible relative Hopf modules maps to the Picard group of $A$, and
the kernel is described as a quotient group
of the group of invertible grouplike elements of the coring $A\ot H$, or as a Harrison
cohomology group. Our methods are based on elementary $K$-theory. The Hilbert 90
Theorem follows as a corollary. The part of the Picard group of the coinvariants
that becomes trivial after base extension embeds in the Harrison cohomology group,
and the image is contained in a well-defined subgroup $E$. It equals $E$ if
$H$ is a cosemisimple Hopf algebra over a field.
\end{abstract}

\maketitle

\section*{Introduction}
Let $l$ be a cyclic Galois field extension of $k$. The Hilbert 90 Theorem
tells us that every cocycle in $Z^1(C_p,l^*)$ is a coboundary. There exist
various generalizations of this result. For example, if we have a Galois
extension $B\to A$ of commutative rings, with Galois group $G$, then the
cohomology group $H^1(G,\GG_m(A))$ is isomorphic to $\Pic(A/B)$, the kernel
of the natural map from the Picard group of $B$ to the Picard group of $A$,
see for example \cite{DI}. Now we can ask the following question: suppose that
$G$ acts on $A$ as a group of isomorphisms. Can we still give an algebraic 
interpretation of $H^1(G,\GG_m(A))$? A second problem is whether there is any relation
between $H^1(G,\GG_m(A))$ and the Picard group of the ring of invariants $B=A^G$.\\
In this note, we will discuss these two problems in a more general situation:
we will assume that $A$ is a commutative $H$-comodule algebra, with $H$
an arbitrary commutative bialgebra over a commutative ring $k$. We then ask
for an algebraic interpretation of the first Harrison cohomology group
$H^1_{\rm Harr}(H,A,\GG_m)$ (with notation as in \cite{C2}). If $H$ is finitely
generated and projective, then this Harrison cohomology group is isomorphic to
a Sweedler cohomology group $Z^1_{\rm Harr}(H,A,\GG_m)$, and if $H=\ZZ G$
with $G$ a finite group, then it reduces to the cohomology group 
$H^1(G,\GG_m(A))$.\\
We proceed as follows: we introduce the relative Picard group $\Pic^H(A)$
as the Grothendieck group of the category of invertible relative Hopf modules.
The forgetful functor to the category of invertible $A$-modules induces
a K-theoretic exact sequence, linking the Picard group of $A$, the relative
Picard group, and the groups of unit elements of $A$ and the coinvariants
$B=A^{{\rm co}H}$; the middle term in the sequence can be computed, and it
is the group of invertible grouplike elements of the coring $A\ot H$.
We show also that these grouplike elements are precisely the Harrison cocycles,
and it follows from the exactness of the sequence that the first Harrison
cohomology group is the kernel of the map $\Pic^H(A)\to \Pic(A)$, answering
our first question.\\
Then we observe that there is a similar exact sequence associated to the induction
functor $\dul{\Pic}(B)\to \dul{\Pic}(A)$, and that the two exact sequences
fit into a commutative diagram. If $A$ is a faithfully flat Hopf Galois extension
of $B$, then the categories of $B$-modules and relative Hopf modules are
equivalent, hence $\Pic(B)\cong\Pic^H(A)$, and we recover Hilbert 90.
In general, we have an injection $\Pic(A/B)\to H^1_{\rm Harr}(H,A,\GG_m)$,
and we can describe a subgroup of $H^1_{\rm Harr}(H,A,\GG_m)$ that contains
the image of $\Pic(A/B)$. The image is precisely this subgroup if 
$H$ is a cosemisimple Hopf algebra over a field $k$.\\
A special situation is the following: let $k$ be an algebraically closed field,
$A$ a finitely generated commutative normal $k$-algebra, and $G$ a connected
algebraic group acting rationally on $A$. Then $A$ is an $H$-comodule algebra,
with $H$ the affine coordinate ring of $G$. In this case, our exact sequence
was given by Magid in \cite{23}, but apparently the author of \cite{23}
was not aware of the connection to
Harrison cohomology, grouplike elements of corings and the generalized Hilbert 90
Theorem.\\
In \seref{4}, we study the Harrison cocycles (or the grouplike
elements in $A\ot H$) in some particular cases. First we look  at the situation
considered by Magid in \cite{23}, and then it turns out that the grouplike
elements of $G(A\ot H)$ are induced by the grouplike elements of $H$.
In the situation where $A$ is a $\ZZ$-graded commutative $k$-algebra,
the relative Picard group turns out to be the graded Picard group $\Pic_g(A)$
studied by the first author in \cite{C1}. If $A$ is reduced, then the grouplike
elements of $A\ot H$ can also be described using the grouplikes of $H$, according
to a result in \cite{C1}.

\section{Preliminary results}\selabel{1}
\subsection{The language of corings}\selabel{1.1}
Relative Hopf modules can be viewed as comodules over a coring. This will
be used in the sequel, and this is why we briefly recall some properties of
corings.
Recall that an $A$-coring is a comonoid in the monoidal category ${}_A\Mm_A$
of $A$-bimodules. Thus an $A$-coring $\CCC$ is an $A$-bimodule together with
two $A$-bimodule maps
$$\Delta_\CCC :\ \CCC\to \CCC\otimes_A\CCC \quad
\hbox{and}\quad \varepsilon_\CCC:\ \CCC\to A$$
satisfying the usual coassociativity and counit properties. We refer to
\cite{4,5,6,15,26} for a detailed discussion of corings.
$$G(\CCC)=\{X\in \CCC~|~\Delta_\CCC(X)=X\ot_A X~{\rm and}~\varepsilon_\CCC(X)=1\}$$
is the set of grouplike elements of $\CCC$. A right $\CCC$-comodule $M$ is a right
$A$-module together with a right $A$-linear map 
$\rho_M :\ M\rightarrow M\otimes_A \CCC$ satisfying
$$(M\otimes_A\varepsilon_\CCC)\circ \rho_M=M, \quad \hbox{and} \quad
(M\otimes_A\Delta_{\CCC})\circ \rho_M=(\rho_M \otimes_A \CCC)\circ \rho_M.$$
A morphism of right $\CCC$-comodules $f :\ M \to N$ is an
$A$-linear map $f$ such that
$$\rho_N \circ f=(f\otimes_A\CCC)\circ \rho_M.$$
$\Mm^\CCC$ will be the category of right $\CCC$-comodules and comodule morphisms.
We have the following interpretation of the grouplike elements of $\CCC$.

\begin{lemma}\lelabel{1.1}
Let $\CCC$ be an $A$-coring. Then there is a bijective correspondence
between $G(\CCC)$ and the set of maps $\rho:\ A\to A\ot_A\CCC=\CCC$ making
$A$ into a right $\CCC$-comodule. The coaction $\rho_X$ corresponding to $X\in G(\CCC)$
is given by
$$\rho_X(a)=Xa.$$
With this notation, $A^X=(A,\rho_X)\cong A^Y= (A,\rho_Y)$ as right $\CCC$-comodules
if and only if there exists an invertible $b\in A$ such that $\rho_Y(b)=Yb=bX$.
\end{lemma}

\begin{proof}
The first part is well-known (and straightforward), see for example \cite{5}.
Let $f:\ A^X\to A^Y$ be a right $\CCC$-colinear isomorphism.
Then $f(a)=ba$ for some $b\in A$, which is
invertible since $f$ is an isomorphism. The fact that $f$ is $\CCC$-colinear
tells us that
$$Yb=\rho_Y(f(1))=(f\ot_A\CCC)(\rho_X(1))=bX.$$
The converse property is obvious.
\end{proof}

\subsection{Relative Hopf modules}\selabel{1.2}
Let $H$ be a bialgebra over a commutative ring $k$, and $A$ a right $H$-comodule
algebra. Throughout this note,
we will assume that $H$ and $A$ are commutative. Then 
$A$ is a commutative algebra
and we have a right $H$-coaction $\rho$ on $A$ such that
$$\rho(ab)=a_{[0]}b_{[0]}\ot a_{[1]}b_{[1]},$$
for all $a,b\in A$. Here we use the Sweedler-Heyneman notation for the coaction $\rho$:
$\rho(a)=a_{[0]}\ot a_{[1]}$, with summation implicitely understood. For the comultiplication
on $H$, we use the notation $\Delta(h)=h_{(1)}\ot h_{(2)}$. \\
A relative Hopf module $M$ is a $k$-module, together with a right $A$-action and a right $C$-coaction
$\rho_M$ such that
$$\rho_M(ma)=m_{[0]}a_{[0]}\ot m_{[1]}a_{[1]},$$
for all $a\in A$ and $m\in M$. The category of relative Hopf modules and $A$-linear $H$-colinear
maps will be denoted by $\Mm^H_A$. The coinvariant submodule $M^{{\rm co}H}$ of $M\in \Mm^H_A$
is defined by
$$M^{{\rm co}H}=\{m\in M~|~\rho_M(m)=m\ot 1\}.$$
$A^{{\rm co}H}=B$ is a $k$-subalgebra of $A$, and $M^{{\rm co}H}$ is a $B$-module. We obtain
a functor $(-)^{{\rm co}H}:\ \Mm^H_A\to \Mm_B$, that has a left adjoint
$T=-\ot_BA:\ \Mm_B\to \Mm^H_A$. The right $H$-coaction on $N\ot_BA$ is $N\ot_B \rho$.
The unit $u$ and counit $c$ of the adjunction are given by the
following formulas, for $N\in\Mm_B$ and $M\in \Mm^H_A$:
$$u_N :\ N \to (N\otimes _{B}A )^{{\rm co}H},~~ u_N(n)=n \otimes 1;$$
$$c_M :\ M^{{\rm co}H}\otimes_{B}A \to M,~~ c_{M}(m\otimes a)=ma.$$
$A$ is called a Hopf algebra extension of $B=A^{{\rm co}H}$ if the canonical map
$${\rm can}:\ A\ot_B A\to A\ot H,~~{\rm can}(a\ot_B b)=ab_{[0]}\ot b_{[1]}$$
is an isomorphism. If $A$ is a faithfully flat Hopf Galois extension, then the adjunction
$(-\ot_B A,(-)^{{\rm co}H})$ is a pair of inverse equivalences.
We refer to \cite{11,24,27} for a detailed discussion of Hopf algebras and relative Hopf modules.\\

 $\CCC=A\ot H$ is a coring, with
structure maps
$$a'(b\ot h)a=a'ba_{[0]}\ot ha_{[1]}$$
$$\Delta_{\CCC}(a\ot h)=(a\ot h_{(1)})\ot_A(1\ot h_{(2)})$$
$$\varepsilon_{\CCC}(a\ot h)=a\varepsilon(h)$$
The category $\Mm^{A\ot H}$ is isomorphic to the category $\Mm_A^H$ of relative Hopf modules;
we refer to \cite{4,6} for full detail. 
Note that
$X=\sum_i a_i \otimes h_i \in G(A\ot H)$ if and only if
\begin{equation}\eqlabel{grouplike1}
\sum_i (a_i \otimes {h_i}_{(1)} \otimes {h_i}_{(2)})=\sum_{i,j} (a_i{a_j}_{[0]} \otimes
h_i{a_j}_{[1]}\otimes h_j)\quad
\hbox{and} \quad \sum a_i\varepsilon(h_i)=1.
\end{equation}
$A\ot H$ is also a commutative algebra, with
multiplication
$$(a\ot h)(b\ot k)=ab\ot hk.$$
The product of two grouplike elements is a grouplike element, and $1_A\ot 1_H$ is grouplike, hence
$G^i(A\ot H)$, the set of invertible grouplike elements, is an abelian group. Also observe that
an invertible grouplike element is precisely a normalized Harrison 1-cocycle (see for example \cite[Sec. 9.2]{C2} for the definition of the Harrison complex).\\
Let $H$ be a finitely generated projective cocommutative Hopf algebra, and $A$ a commutative left $H$-module algebra. Then $H^*$
is a commutative Hopf algebra, and $A$ is a right $H^*$-comodule algebra. 
If $\sum_i a_i\ot f_i\in A\ot H^*$ is an invertible grouplike element (or a normalized Harrison cocycle),
then 
\begin{equation}\eqlabel{Sweedler1}
\phi:\ H\to A,~~\phi(h)=\sum_i a_if_i(h),
\end{equation}
is a normalized Sweedler $1$-cocycle, this means that $\phi(1_H)=1_A$, and the cocycle condition
\begin{equation}\eqlabel{Sweedler2}
\phi(hh')=\sum_i(h_{(1)}.(\phi(h')))\phi(h_{(2)}),
\end{equation}
is satisfied. This gives a bijective correspondence between Harrison and Sweedler cocycles,
see \cite[Prop. 9.2.3]{C2}. For the definition of the Sweedler complex, see \cite{S2} or
\cite[Sec. 9.1]{C2}. In the case where $H=kG$, with $G$ a finite group, Sweedler cohomology reduces to
group cohomology.

\subsection{Elementary algebraic K-theory}\selabel{1.3}
Let $(\Cc,\ot, I)$ and $(\Dd,\ot, J)$ be skeletally small symmetric monoidal categories,
and let $F:\ \Cc\to \Dd$ be a cofinal, strong monoidal functor. Then we can consider
the Grothendieck and Whitehead groups of $\Cc$ and $\Dd$, and we have an exact sequence
connecting them (see for example \cite[Ch. VII]{Bass}):
\begin{equation}\eqlabel{K1}
K_1\Cc\rTo^{K_1F} K_1\Dd\rTo^{d}K_1\ul{\phi F}\rTo^{g}
K_0\Cc\rTo^{K_0F} K_0\Dd.
\end{equation}
$C\in \Cc$ is called invertible if there exists $C'\in \Cc$ such that $C\ot C'\cong I$. If all elements
of $\Cc$ and $\Dd$ are invertible, then the description of the five groups in \equref{K1} and
the connecting maps simplifies (see \cite[App. C]{C2}). $K_0\Cc$ is the group of isomorphism
classes of objects in $\Cc$ and $K_1\Cc\cong \Aut_{\Cc}(I)$ (which is then an abelian group).
Let $\ul{\Psi F}$ be the following category: objects are couples $(C,\alpha)$, with
$C\in \Cc$ and $\alpha:\ F(C)\to J$ an isomorphism in $\Dd$.  A morphism between
$(C,\alpha)$ and $(C',\alpha')$ is an isomorphism $f:\ C\to C'$ in $\Cc$ such that
$\alpha'= F(f)\circ \alpha$. $\ul{\Psi F}$ is monoidal, every object is invertible and
$K_1\ul{\phi F}\rTo^{g}\cong K_0\ul{\Psi F}$. The maps $d$ and $g$ are given as follows:
$d(\alpha)=[(I,\alpha)]$ and $g[(C,\alpha)]=[C]$.\\
A typical example is the following: For a commutative ring $A$, 
let $\dul{\Pic}(A)$ be the category of invertible $A$-modules.
If $i:\ B\to A$ is a morphism of commutative rings, then we have cofinal strongly monoidal functor
$$G=-\ot_B A:\  \dul{\Pic}(B)\to \dul{\Pic}(A),$$
and \equref{K1} takes the form
\begin{equation}\eqlabel{K2}
1\to \GG_m(B)\to\GG_m(A)\rTo^{d'}K_1\ul{\phi G}\rTo^{g'} \Pic(B)\to \Pic(A).
\end{equation}

\section{The relative Picard group}\selabel{2}
If $M,N\in \Mm_A^H$, then $M\ot_A N\in \Mm_A^H$, with right $H$-coaction
$$\rho_{M\ot_AN}(m\ot_An)=m_{[0]}\ot_A n_{[0]}\ot m_{[1]}n_{[1]}.$$
So we have a symmetric monoidal category $(\Mm_A^H,\ot_A,A)$. Let $\dul{\Pic}^H(A)$
be the full subcategory consisting of invertible objects.
$\Pic^H(A)=K_0\dul{\Pic}^H(A)$, the group of isomorphism classes of relative Hopf modules,
will be called the relative Picard group of $A$ and $H$. The isomorphism class in $\Pic^H(A)$
represented by an invertible relative Hopf module $M$ will be denoted by $\{M\}$.
This new invariant fits into an exact sequence:

\begin{proposition}\prlabel{2.1}
We have an exact sequence
\begin{equation}\eqlabel{2.1.1}
1\to \GG_m(B)\to\GG_m(A)\rTo^{d}G^i(A\ot H)\rTo^{g} \Pic^H(A)\to \Pic(A).
\end{equation}
\end{proposition}

\begin{proof}
This result can be proved in two ways: a first possibility is to show that \equref{2.1.1}
is precisely the exact sequence \equref{K1}, associated to the functor
$\dul{\Pic}^H(A)\to \dul{\Pic}(A)$ forgetting the $H$-coaction. Let us present an
easy direct proof.\\
The map $\GG_m(B)\to\GG_m(A)$ is the natural inclusion. Take $a\in A$
invertible, and let $d(a)=X=a^{-1}a_{[0]}\ot a_{[1]}$. $X$ is grouplike, since
$a^{-1}a_{[0]}\varepsilon(a_{[1]})=1$, and
\begin{eqnarray*}
X\ot_A X&=&
(a^{-1}a_{[0]}\ot a_{[1]})\ot_A(b^{-1}b_{[0]}\ot b_{[1]})\\
&=&a^{-1}a_{[0]}(b^{-1})_{[0]}b_{[0]}\ot a_{[1]}(b^{-1})_{[1]}b_{[1]}\ot b_{[2]}\\
&=&a^{-1}b_{[0]}\ot b_{[1]}\ot b_{[2]}\\
&=& (a^{-1}b_{[0]}\ot b_{[1]})\ot_A(1\ot b_{[2]})=\Delta(X),
\end{eqnarray*}
where we identified $(A\ot H)\ot_A(A\ot H)=A\ot H\ot H$ and we wrote $a=b$.
The inverse of $X$ is $X^{-1}=a(a^{-1})_{[0]}\ot (a^{-1})_{[1]}$, so $X\in G^i(A\ot H)$.\\
If $d(a)=a^{-1}a_{[0]}\ot a_{[1]}=1_A\ot 1_H$, then $a_{[0]}\ot a_{[1]}=a\ot 1_H$, 
so $a\in B$, and the
sequence is exact at $\GG_m(A)$.\\
For $X\in G^i(A\ot H)$, let $g(X)=A^X$, with notation as in \leref{1.1}. 
$g$ is multiplicative: take $X=\sum_i a_i\ot h_i$ and $Y=\sum_j b_j\ot k_j$
in $G^i(A\ot H)$, then $A^X\ot_AA^Y=A$ as an $A$-bimodule, with comultiplication
given by $\rho_{A^X\ot_AA^Y}(1)=\sum_{i,j} a_i\ot_A b_j\ot h_ik_j=XY$, as needed.\\
If $g(X)=\{A\}$ in $\Pic^H(A)$, then there
exists an $H$-colinear $A$-linear isomorphism
$f:\ A^X\to A$. Then $f(1)=a$ is invertible in $A$, and, since $f$ is $H$-colinear,
$a_{[0]}\ot a_{[1]}=\rho(a)=(f\ot H)(X)=aX$, so $X=a^{-1}a_{[0]}\ot a_{[1]}=d(a)$, and the sequence
is also exact at $G^i(A\ot H)$.\\
The exactness of the sequence at $\Pic^H(A)$ follows from \leref{1.1}.
\end{proof}

\begin{remark}\relabel{2.2}
Let $H=k\ZZ$, and let $A$ be a commutative $\ZZ$-graded $k$-algebra.
Then $\Pic^H(A)=\Pic_g(A)$, the graded Picard group of $A$, as introduced
in \cite{C1}, see also \cite{CVO}. The exact sequence \equref{2.1.1} reduces
to the exact sequence in \cite[Prop. 2.1]{C1}.
\end{remark}

The map $d:\ \GG_m(A)\rTo^{d}G^i(A\ot H)$ is precisely the map $\GG_m(A)\to \GG_m(A\ot H)$
in the Harrison complex. From \prref{2.1}, we therefore obtain immediately:

\begin{corollary}\colabel{2.3}
With $H$ and $A$ as in \prref{2.1}, we have an isomorphism of abelian groups
$$\Pic^H(A)\cong H^1_{\rm Harr}(H,A,\GG_m).$$
\end{corollary}

This is the promised algebraic interpretation of the first Harrison cohomology group.
Note that there are no flatness or projectivity assumptions on $H$ or $A$. We have
Hilbert 90 as an easy consequence. 

\begin{corollary}\colabel{2.4} {\bf (Hilbert 90)}
Let $H,A,B$ be as in \prref{2.1}. If $A$ is a faithfully flat $H$-Galois extension of $B$,
then we have an isomorphism of abelian groups
$$\Pic(A/B)\cong H^1_{\rm Harr}(H,A,\GG_m).$$
\end{corollary}

\begin{proof}
From the fact that the monoidal categories $\Mm_B$ and $\Mm^H_A$ are equivalent, it
follows that $\Pic(B)\cong\Pic^H(A)$.
\end{proof}

Take the exact sequences \equref{K2} and \equref{2.1.1},
and observe that they fit into a commutative diagram
$$\begin{diagram}
1&\rTo^{} &\GG_m(B)&\rTo^{}&\GG_m(A)&\rTo^{d'}&K_1\ul{\phi G}&\rTo^{g'}& \Pic(B)&\rTo^{}& \Pic(A)\\
&&\dTo_{=}&&\dTo_{=}&&&&\dTo_{j}&&\dTo_{=}\\
1&\rTo^{} &\GG_m(B)&\rTo^{}&\GG_m(A)&\rTo^{d}&G^i(A\ot H)&\rTo^{g}& \Pic^H(A)&\rTo^{}& \Pic(A)
\end{diagram}$$
The map $j$ maps $[N]\in \Pic(B)$ to $\{N\ot_B A\}\in \Pic^H(A)$.
Using the Five Lemma, we find a map $i:\ K_1\ul{\phi G}\to G^i(A\ot H)$.

\begin{lemma}\lelabel{2.5}
With notation as above, the maps $i$ and $j$ are injective.
\end{lemma}

\begin{proof}
From the fact that $u$ is a natural transformation between additive endofunctors
of the category of $B$-modules, and since $u_B$ is
an isomorphism, it follows that $u_N:\ N\to
(N\ot_B A)^{{\rm co}H}$ is an isomorphism
if $N$ is finitely generated and projective as a $B$-module. So if $N\ot_B A\cong A$, then
$N\cong (N\ot_B A)^{{\rm co}H}\cong A^{{\rm co}H}=B$, and $j$ is injective.
The injectivity of $i$ then follows from an easy diagram chasing argument.
\end{proof}

Our next aim is to characterize the image of $i$. This will be the topic of \seref{3};
it will turn out that we obtain nice results in the case where $H$ is cosemisimple.

\section{Coinvariantly generated relative Hopf modules}\selabel{3}
Some of our results will be more specific if we assume
that $H$ is a cosemisimple Hopf algebra over a field $k$. Recall that $H$ is cosemisimple if
there exists a left integral
$\phi$ on $H^*$ such that $\phi(1)=1$ (see e.g. \cite{26}). In this case, the coinvariants
functor $(-)^{{\rm co}H}:\ \Mm_A^H\to \Mm_B$ is exact, see \cite[Lemma 2.4.3]{24}.

A relative Hopf module $M$ is called {\sl coinvariantly generated} if $c_M$
is surjective, or, equivalently, if $M=M^{{\rm co}H}A$. If $M$ is coinvariantly
generated, and finitely generated as an $A$-module, then we can find a finite
set $\{m_1,\cdots,m_n\}\in M^{{\rm co}H}$ that generates $M$.

It follows immediately
from the properties of adjoint functors that $N\ot_BA$ is coinvariantly generated,
for every  $N\in \Mm_B$; in particular, $A$ is coinvariantly generated.
We also have the following:

\begin{lemma}\lelabel{3.1} 
Let $M\in \Mm^H_A$ and $N\in \Mm_B$. If $M$ is an
epimorphic image of $N\otimes_BA$ in
$\Mm^H_A$, then $M^{{\rm co}H}=0$ implies that $M=0$.
\end{lemma}

\begin{proof} If $M^{{\rm co}H}=0$, then $\Hom_A^H(N\otimes_B A, M)=\Hom_B(N, M^{{\rm co}H})=0$. But $\Hom_A^H(N\otimes_B A, M)$ contains the
epimorphism of relative Hopf modules  $N\otimes_B A\rightarrow M$, so $M=0$.
\end{proof}

If $N$ is an epimorphic image of $M$ in $\Mm^H_A$, and if $M$ is
coinvariantly generated, then $N$ is also coinvariantly generated.

\begin{lemma}\lelabel{3.2}
Assume that
$H$ is a cosemisimple Hopf algebra over a field $k$. If $N\in \Mm_B$ is
projective, then $N\ot_B A$ is projective in $\Mm^H_A$.
\end{lemma}

\begin{proof} 
See \cite[Prop. 2.5]{9}.
\end{proof}

\begin{lemma}\lelabel{3.3} Let $k$ be a field.
\begin{enumerate}
\item The forgetful functor $\Mm^H_A\to\Mm_A$ preserves projectives;
\item if $H$ is cosemisimple, then the forgetful functor also reflects
projectivity of finitely generated modules.
\end{enumerate}
\end{lemma}

\begin{proof}
(1) Take $M\in \Mm^H_A$ projective, and consider the epimorphism
$p:\ M\ot A\to M$, $p(m\ot a)=ma$ in $\Mm^H_A$. The exact sequence
$$0\to \Ker p\to M\ot A\rTo^{p}M\to 0$$
splits in $\Mm^H_A$, since $M$ is a projective object, and a fortiori
in $\Mm_A$. Hence $M$ is a direct factor of $M\ot A$, which is a projective
 $A$-module, so $M$ is also a projective  $A$-module.\\
(2) Let $M,N$ be a relative Hopf modules, and assume that $M$ is finitely generated
and projective in $\Mm_A$. According to \cite[Prop. 4.2]{9},
 $\Hom_A(M,N)\in \Mm_A^H$, and
it is easy to show that $\Hom_A(M,N)^{{\rm co}H}=\Hom_A^H(M,N)$. It follows that
the functor $\Hom_A^H(M,-):\ \Mm_A^H\to \Mm$ is exact, since it is the composition of
the exact functors 
$\Hom_A(M,-):\ \Mm^H_A\to \Mm^H$ ($M\in \Mm_A$ is projective)
and $(-)^{{\rm co}H}:\ \Mm^H\to \Mm$ ($H$ is cosemisimple).
\end{proof}

\begin{lemma}\lelabel{3.4}
Let $H$ be a cosemisimple Hopf algebra over a field $k$, and take $P,Q\in \Mm^H_A$ finitely generated as $A$-modules. Assume that $Q$ is a projective object of $\Mm^H_A$.
Then every epimorphism $f:\ P\to Q$ in $\Mm^H_A$ has a right inverse
in $\Mm^H_A$.
\end{lemma}

\begin{proof}
It is clear that $\Hom_A(Q,P)$ and $\Hom_A(Q,Q)$ are right $H$-comodules,
and the map $f^*=\Hom_A(Q,f):\ \Hom_A(Q,P)\to \Hom_A(Q,Q)$ is right
$H$-colinear. It follows from \leref{3.3} that $Q$ is projective as an 
$A$-module, so $f^*$ is surjective. Since $f^*$ is $H$-colinear, it restricts
to a surjection 
$$\Hom_A^H(Q,P)=\Hom_A(Q,P)^{{\rm co}H}\to \Hom_A^H(Q,Q)=\Hom_A(Q,Q)^{{\rm co}H}.$$
Take a preimage $g\in \Hom_A^H(Q,P)$ of the identity map ${\rm id}_Q$ on $Q$. Then
$f\circ g={\rm id}_Q$, and the result follows.
\end{proof}

For $M\in \Mm_A$, we will denote the dual module by $M^*=\Hom_A(M,A)$.

\begin{proposition}\prlabel{3.5}
Let $H$ be  cosemisimple, and assume that $P\in \Mm^H_A$ is 
coinvariantly generated and finitely generated projective as an $A$-module.
Then
\begin{enumerate}
\item $P^{{\rm co}H}$ is a finitely generated projective $B$-module;
\item $P^*$ is coinvariantly generated;
\item the map $c_P$ is an isomorphism in $\Mm^H_A$.
\end{enumerate}
\end{proposition}

\begin{proof} (1) As we have seen, 
there exist $p_1, p_2,\cdots, p_n \in P^{{\rm co}H}$
such that $P=\sum_i p_iA$. Set $F=A^n$ and
let $f:\ F\ \to P$ be the
$A$-linear map given  by $f(a_1, a_2, ... , a_n)=\sum_i p_ia_i$. Then 
$F\in \Mm^H_A$
and $f$ is an epimorphism in $\Mm^H_A$. By \leref{3.4}, there
exists a monomorphism $g\in \Hom_A(P ,
F)$ such that $f\circ g={\rm id}_P$. The restriction of $g$ to $P^{{\rm co}H}$ is then
a $B$-linear
right inverse of the restriction of $f$ to $F^{{\rm co}H}$, and $F^{{\rm co}H}=B^{n}$,
and we obtain (1).

(2) The map $g^*=\Hom_A(g , A):\ F^*\to P^*$ is surjective and
$H$-colinear. The fact that $F^*$ is coinvariantly generated then implies that
$P^*$ is also coinvariantly generated.

(3) Consider the natural transformation $t:\ (-)^{{\rm co}H}\ot_BA\to (-)$
given by
$$t_P:\ P^{{\rm co}H} \ot_BA\rightarrow P,~~t_P(p\otimes a)= pa.$$
The map $t_A$ is an isomorphism, so
$t_F$ is an isomorphism by
additivity. It follows that $t_P$ is an isomorphism, since $F=P\oplus \Ker
f$ as $H$-comodules.
\end{proof}

Let $X=\sum_i a_i\ot h_i\in G(A\ot H)$, and write
$$A_X=\{a\in A~|~\rho(a)=aX=\sum_i aa_i\ot h_i\},$$
and
$$A_X^i=\{a\in A_X~|~a~{\rm is~invertible}\}.$$
Observe that
$$\im(d)=\{X\in G^i(A\ot H)~|~A_X^i\neq \emptyset\}$$
and
$$A_{1\ot 1}=A^{{\rm co}H}.$$
Furthermore, $A_XA_Y\subset A_{XY}$: take $a\in A_X$ and $b\in A_Y$, then
$\rho(a)=aX=\sum_i aa_i\ot h_i$, $\rho(b)=bY=\sum_j bb_j\ot k_j$ and
$\rho(ab)=a_{[0]}b_{[0]}\ot a_{[1]}b_{[1]}=\sum_{i,j} aa_ibb_j\ot h_ik_j=abXY$. Also
$A_X^i\cap A_Y^i=\emptyset$ if $X\neq Y$.

\begin{lemma}\lelabel{3.6}
The set
$$E=\{X\in G^i(A\ot H)~|~AA_X=A~{\rm and}~AA_{X^{-1}}=A\}$$
is a subgroup of $G^i(A\ot H)$ containing $\im(d)$.
\end{lemma}

\begin{proof}
If $X\in \im(d)$, then there exists an invertible $a\in A_X$, and then
$AA_X=A$. since $X^{-1}\in \im (d)$, we also have $AA_{X^{-1}}=A$, hence
$X\in E$. It is clear that $1\ot 1\in E$. If $X,Y\in E$, then
$AA_{XY}\supset AA_XA_Y=AA_Y=A$, and, in a similar way,
$AA_{(XY)^{-1}}=A$, hence $XY\in E$. Finally, if $X\in E$, then obviously
$X^{-1}\in E$.
\end{proof}

\begin{proposition}\prlabel{3.7}
Consider the injective map $j:\ \Pic(B)\to \Pic^H(A)$.
If $H$ is a cosemisimple Hopf algebra over a field $k$, then
$$\im(j)=\{\{M\}\in \Pic^H(A)~|~M~\hbox{is coinvariantly generated}\}.$$
\end{proposition}

\begin{proof}
$M\ot_B A$ is coinvariantly generated, so $\im (j)$ is contained in the
desired set. If $H$ is cosemisimple, and $\{N\}\in \Pic^H(A)$, with
$N$ coinvariantly generated, then $N=(N^{{\rm co}H})\ot_BA\in \im (j)$,
by \prref{3.5}(3).
\end{proof}

\begin{lemma}\lelabel{3.8}
Take $X\in G^i(A\ot H)$. Then $A^X$ is coinvariantly generated if and only if
$AA_{X^{-1}}=A$. If $H$ is cosemisimple, then this is also equivalent to
$X\in E$.
\end{lemma}

\begin{proof}
The first statement follows from the fact that $(A^X)^{{\rm co}H}=A_{X^{-1}}$.
Indeed, $a\in (A^X)^{{\rm co}H}$ if and only if $\rho_X(a)=Xa=a\ot 1$, if and only if
$\rho(a)=(1\ot 1)a=X^{-1}(a\ot 1)=aX^{-1}$, which means that $a\in A_{X^{-1}}$.\\
Let $H$ be cosemisimple. Note that $(A^X)^*\cong A^{X^{-1}}$ as relative Hopf
modules. If $A^X$ is coinvariantly generated, then so is $A^{X^{-1}}$, by
\prref{3.5}, and then $X\in E$.
\end{proof}

Now we are able to prove the main result of this Section.

\begin{theorem}\thlabel{3.9}
Consider the monomorphism 
$i:\ K_1\ul{\phi G}\to G^i(A\ot H)$ introduced in \leref{2.5}.\\
Then $\im(i)\subset E$ and $\im(i)= E$ if $H$ is a cosemisimple Hopf algebra over
a field $k$. In this situation, $\Pic(A/B)\cong E$.
\end{theorem}

\begin{proof}
Take $[(M,\alpha)] \in 
K_0\ul{\psi G}$, and let $i[(M,\alpha)]=X\in G^i(A\ot H)$. Then
$\{A^X\}=j(g'[(M,\alpha)])=j([M])=\{M\ot_BA\}$, hence $A^X$ is coinvariantly
generated and $AA_{X^{-1}}=A$, by \leref{3.8}. In a similar way,
$i([(M,\alpha)]^{-1})=X^{-1}$, and $A^{X^{-1}}\cong M^*\ot_B A$ is coinvariantly
generated, so $AA_X=A$, again by \leref{3.8}. This proves that $X\in E$.\\
Assume now that $H$ is cosemisimple, and take $X\in E$. It follows from
\leref{3.8} that $A^X$ is coinvariantly generated, and from \prref{3.7}
that $A^X=M\ot_B A$ for some $M\in \dul{\Pic}(B)$. Since the image of $M$
in $\Pic(A)$ is trivial, $[M]=g'[(M,\alpha)]$ for some $(M,\alpha)\in \Cc$.
Write $i[(M,\alpha)]=Y$. Then $X=Yd(a)$, for some $a\in \GG_m(A)$. Consider the map $\alpha':\ M\ot_BA\to A$, $\alpha'(m\ot b)=a^{-1}\alpha(m\ot b)$.
Then $i[(M,\alpha')]=X$.
\end{proof}

\section{On the grouplike elements}\selabel{4}
We have an injective map $i:\ G(H)\to G(A\ot H)$, $i(g)=1_A\ot g$. Everything
simplifies if $i$ is an isomorphism. We discuss two situations in which this
is (almost) the case.

Recall that a commutative algebra which is
an integral domain is called normal if it is integrally closed in its field
of fractions.

\begin{proposition}\prlabel{4.11}
 Let $k$ be an algebraically closed field, $A$ a finitely generated commutative
normal $k$-algebra and $G$ a
connected algebraic group acting rationally on $A$. Let $H$ be the affine
coordinate ring of $G$,
and $\chi(G)$ be the group of characters of $G$. Then
$$G(A\otimes H)=\{1 \otimes \phi~|~\phi \in  G(H)=\chi(G)\}.$$
\end{proposition}

\begin{proof} Let $x=\sum_i a_i \otimes f_i \in G(A \otimes H)$. Then we have
\begin{equation}\eqlabel{4.11.1}
\sum_i (a_i \otimes {f_i}_{(1)} \otimes {f_i}_{(2)})=\sum_{i,j} (a_i{a_j}_{[0]} \otimes
(f_i*{a_j}_{[1]})\otimes f_j)
\end{equation}
and $\sum a_i\varepsilon(f_i)=1$.
The map
$$\alpha:\ A\ot H\to\Hom(kG,A),~~\alpha(a\ot f)(g)=af(g)$$
is injective. Let $\phi=\alpha(x)$. 
Using \equref{4.11.1},
we compute for all $g,g'\in G$ that
\begin{eqnarray*}
&&\hspace*{-15mm}
\phi(gg')=\sum_i a_if_i(gg')=\sum_i a_i{f_i}_{(1)}(g){f_i}_{(2)}(g')\\
&=&\sum_{i,j} (a_i{a_j}_{[0]}((f_i*{a_j}_{[1]})(g))f_j(g')=
\sum_{i,j} a_i{a_j}_{[0]}f_i(g){a_j}_{[1]}(g)f_j(g')\\
&=&\sum_{i,j} (g.a_j)f_j(g')a_if_i(g)=\sum_{i,j}
g.(a_jf_j(g'))a_if_i(g)\\
&=&(g.(\phi(g')))\phi(g).
\end{eqnarray*}
From the second equality, we have $1=\sum_i a_if_i(1_G)=\phi(1_G)$. For every $g\in G$,
$\phi(g)$ is invertible in
$A$, with inverse $g.(\phi(g^{-1}))$. By the proof of
\cite[Prop. 1b, p. 46]{21}, $\phi(g)\in k$ for every $g\in G$, so $\phi \in\chi(G)$. 
Now $\chi(G)=G(H)\subset H$ (see \cite[p. 25]{19}), so it follows in particular that $\phi\in H$.
For all $g\in G$ we now have that
$$\alpha(1\ot \phi)(g)=\phi(g)=\sum_i a_if_i(g)=\alpha(x)(g),$$
hence $x=1\ot \phi$, by the injectivity of $\alpha$.
\end{proof}

Now consider the situation from \reref{2.2}: $H=k\ZZ\cong k[X,X^{-1}]$, and $A$ is a commutative
$\ZZ$-graded algebra. In this situation $A\ot H=A\ot k[X,X^{-1}]$. Grouplike elements
in $A\ot H$ can be constructed as follows. Let $1=e_1+\cdots+e_n$ with the $e_i$
orthogonal idempotents, and take $d_1,\cdots, d_n\in \ZZ$. Then 
$\sum_{i=1}^n e_i\ot X^{d_i}$ is a grouplike element in $A\ot k[X,X^{-1}]$.
In this way, we have an embedding of $\Cc(\Spec(A),\ZZ)$, the continuous functions
from $\Spec(A)$ (with the Zariski topology) to $\ZZ$ (with the discrete topology),
into $G(A\ot k[X,X^{-1}])$. The first author was amazed to see that one of his first results,
\cite[Theorem 2.3]{C1},  can be restated in such a way that it becomes a result about
corings. Recall that a commutative ring is called reduced if it has no nontrivial
nilpotents.

\begin{proposition}\prlabel{gr}
Let $A$ be a reduced $\ZZ$-graded commutative $k$-algebra. Then
the map $\Cc(\Spec(A),\ZZ)\to G(A\ot k[X,X^{-1}])$ is a bijection.
\end{proposition}

\begin{example}\exlabel{gr2} (cf. \cite[Example 2.6]{C1})
\prref{gr} does not hold if $A$ contains nilpotent elements; this is related
to the fact that there exist non-homogeneous units in this situation.
Let $A=k[x]$, with $x^2=0$, and put a $\ZZ$-grading on $A$ by taking $\deg (x)=1$.
Then $1+ax\in \GG_m(A)$, and $d(1+ax)=(1-ax)\ot 1+ax\ot X$ is a grouplike
element in $G(A\ot k[X,X^{-1}])$ which is not in the  image of 
$\Cc(\Spec(A),\ZZ)$.
\end{example}


\begin{thebibliography}{99}
\bibitem{Bass}
H. Bass, ``Algebraic K-Theory", Benjamin, New York, 1968.

\bibitem{4} 
T. Brzezi\'nski, The structure of corings. Induction functors, Maschke-type
theorem, and Frobenius and Galois properties,
{\sl Algebr. Representat. Theory} {\bf 5} (2002), 389--410.

\bibitem{5} 
T. Brzezi\'nski, The structure of corings with a grouplike
element, {\sl Banach Center Publ.} {\bf 61} (2003), 21--35.

\bibitem{6} 
T. Brzezi\'nski and R. Wisbauer, ``Corings and comodules", 
{\sl London Math. Soc. Lect. Note Ser.} {\bf 309},
  Cambridge University Press, Cambridge, 2003.
  
\bibitem{C1}
S. Caenepeel, A cohomological interpretation of the graded Brauer group I,
{\sl Comm. Algebra} {\bf 11} (1983), 2129--2149.

\bibitem{C2}
S. Caenepeel, Brauer groups, Hopf algebras and Galois theory, 
{\sl K-Monographs Math.} {\bf 4}, Kluwer Academic Publishers, Dordrecht, 1998.

\bibitem{9} 
S. Caenepeel and T. Gu\'ed\'enon, Projectivity of a
relative Hopf module over the
subring of coinvariants, in ``Hopf algebras", J.
Bergen, S. Catoiu and W. Chin (eds.), {\sl Lect.
Notes Pure Appl. Math.} {\bf 237}, Marcel Dekker, New York, 2004, 97--108.

\bibitem{CVO}
S. Caenepeel and F. Van Oystaeyen, Brauer groups and the cohomology of graded rings,
{\sl Monographs and Textbooks in Pure and Appl. Math.} 
{\bf 121}, Marcel Dekker, New York, 1988.

\bibitem{11} 
S. D\v asc\v alescu, C. N\v ast\v asescu and \c S. Raianu,
``Hopf algebras: an Introduction'', {\sl Monographs Textbooks in Pure
Appl. Math.} {\bf 235}, Marcel Dekker, New York, 2001.

\bibitem{DI}
F. DeMeyer, E. Ingraham, Separable algebras over commutative rings, 
{\sl Lecture Notes in Math.} {\bf 181}, Springer Verlag, Berlin, 1971.

\bibitem{15} 
L. El Kaoutit, J. G\'omez-Torrecillas and F.J.
Lobillo, Semisimple corings, {\sl Algebra Colloquium}, to appear.

\bibitem{19} 
J. C. Jantzen, ``Representations of algebraic groups'',
{\sl Pure Appl. Math.} {\bf 131},  Academic Press, Boston, 1987.

\bibitem{21} 
A. Magid, Finite generation of class groups of rings of
invariants, {\sl Proc. Amer. Math.
Soc.} {\bf 60} (1976), 45--48.

\bibitem{23} 
A. Magid, Picard groups of rings of invariants, {\sl J. Pure Appl.
Algebra} {\bf 17} (1980), 305--311.

\bibitem{24} 
S. Montgomery, ``Hopf algebras and their actions on rings", American
Mathematical Society, Providence, 1993.

\bibitem{S2}
M. E. Sweedler, Cohomology of algebras over Hopf algebras,
{\sl Trans. Amer. Math. Soc.} {\bf 133} (1968), 205-239.

\bibitem{27} 
M. E. Sweedler, ``Hopf algebras'', Benjamin, New York, 1969.

\bibitem{26}
M. E. Sweedler, The predual Theorem to the Jacobson-Bourbaki Theorem,
{\sl Trans. Amer. Math. Soc.} {\bf 213} (1975), 391--406.


\end{thebibliography}
\end{document}